\documentclass[english,11pt]{article}
\usepackage[latin1]{inputenc}
\usepackage{amsmath,amsthm,amssymb,babel}

\newtheorem{teo}{Theorem}

\newtheorem{prop}{Proposition}

\newtheorem{lemma}{Lemma}

\def\proof{{\it Proof.}\ }
\def\endproof{\hfill $\Box$\par\vskip3mm}

\def\eq#1{(\ref{#1})}

\def\neweq#1{\begin{equation}\label{#1}}
\def\endeq{\end{equation}}

\def\phi{\varphi}
\def\RR{{\mathbb R} }

\def\di{\displaystyle}

\date{}

\title{\sc  Neumann problems associated to nonhomogeneous differential operators in Orlicz--Sobolev spaces\thanks{
Correspondence address: Vicen\c{t}iu R\u{a}dulescu, Department of
Mathematics, University of Craiova,  200585 Craiova, Romania. E-mail:
{\tt
vicentiu.radulescu@math.cnrs.fr}}}
\author{\sc Mihai Mih\u ailescu$\,^a$ and Vicen\c{t}iu R\u{a}dulescu$\,^{a,b}$\\
\small
$^a\,$Department of Mathematics, University of Craiova,  200585 Craiova,
Romania\\
\small $^b\,$Institute of Mathematics ``Simion Stoilow" of the Romanian Academy,\\
\small P.O. Box 1-764, 014700 Bucharest, Romania\\
\small
E-mail addresses: {\tt mmihailes@yahoo.com}\qquad {\tt
vicentiu.radulescu@math.cnrs.fr}}

\begin{document}
\baselineskip16pt
\maketitle

\noindent{\small{\sc Abstract}. We study a nonlinear Neumann
boundary value problem associated to a nonhomogeneous differential
operator. Taking into account the competition between the
nonlinearity and the bifurcation parameter, we
establish sufficient conditions for the existence of nontrivial solutions in a related Orlicz--Sobolev space.
}

\section{Introduction and preliminary results}
This paper is motivated by phenomena which are described by
nonhomogeneous Neumann problems of the type
\begin{equation}\label{1}
\left\{\begin{array}{lll}
-{\rm div}(a(x,|\nabla u(x)|)\nabla u(x))+a(x,|u(x)|)u(x)=\lambda\; g(x,u(x)), &\mbox{for}&
x\in\Omega\\
\di\frac{\partial u}{\partial\nu}(x)=0, &\mbox{for}& x\in\partial\Omega\,,
\end{array}\right.
\end{equation}
where $\Omega$ is a bounded domain in $\RR^N$ ($N\geq 3$) with
smooth boundary $\partial\Omega$ and $\nu$ is the outward unit
normal to $\partial\Omega$. In \eq{1} there are also involved the
functions $a(x,t)$,
$g(x,t):\overline\Omega\times\RR\rightarrow\RR$ which will be
specified later and the constant $\lambda>0$.

In the particular case when in \eq{1} we have $a(x,t)=t^{p(x)-2}$,
with $p(x)$ a continuous function on $\overline\Omega$, we deal
with problems involving variable growth conditions. The study of
such problems has been stimulated by  recent advances in
elasticity (see \cite{Z1,Z2}), fluid dynamics (see
\cite{R,RR,D,hal}), calculus of variations and differential
equations with $p(x)$-growth conditions (see
\cite{AM,marcel,mihradcras,RoyalSoc,mihradproc,mihradjmaa,mihradmanuscr,Z1,Z2}).

Another recent application which uses operators as those described
above can be found in the framework of image processing. In that
context we refer to the study of Chen, Levine and Rao \cite{CLR}.
In \cite{CLR} the authors study a functional with variable
exponent, $1<p(x)<2$, which provides a model for image
restoration.  The diffusion resulting from the proposed model is a
combination of Gaussian smoothing and regularization based on
Total Variation. More exactly, the following adaptive model was
proposed
\begin{equation}\label{stea}
\min_{I=u+v,\;u\in{\rm BV}\cap L^2(\Omega)}\int_\Omega\phi(x,\nabla u)\;dx+\lambda\cdot\|u\|^2_{L^2(\Omega)}\,,
\end{equation}
where $\Omega\subset\RR^2$ is an open domain,
$$\phi(x,r)=
\left\{\begin{array}{lll}
\frac{1}{p(x)}|r|^{p(x)}, &\mbox{for}& |r|\leq\beta\\
|r|-\frac{\beta\cdot p(x)-\beta^{p(x)}}{p(x)}, &\mbox{for}& |r|>\beta\,,
\end{array}\right.$$
where $\beta>0$ is fixed and $1<\alpha\leq p(x)\leq 2$. The function $p(x)$ involved here depends on the location
$x$ in the model. For instance it can be used
$$p(x)=1+\frac{1}{1+k|\nabla G_\sigma\ast I|^2}\,,$$
where $G_\sigma(x)=\frac{1}{\sigma}\exp(-|x|^2/(4\sigma^2))$ is the Gaussian filter and $k>0$ and $\sigma>0$ are 
fixed
parameters (according to the notation in \cite{CLR}). For problem \eq{stea} Chen, Levine and Rao establish the 
existence
and uniqueness of the solution and the long-time behavior of the associated flow of the proposed model. The 
effectiveness
of the model in image restoration is illustrated by some experimental results included in the paper.

We point out that the model proposed by Chen, Levine and Rao in
problem \eq{stea} is linked with the energy which can be
associated with problem \eq{1} by taking $\phi(x,\nabla
u)=a(x,|\nabla u|)\nabla u$. Furthermore, the operators which will
be involved in problem \eq{1} can be more general than those
presented in the above quoted model. That fact is due to the
replacement of $|t|^{p(x)-2}t$ by more general functions
$\phi(x,t)=a(x,|t|)t$. Such functions will demand some new setting
spaces for the associated energy, the {\it generalized
Orlicz-Sobolev spaces} $L^\Phi(\Omega)$, where
$\Phi(x,t)=\int_0^t\phi(x,s)\;ds$. Such spaces originated with
Nakano \cite{nak} and were developed by Musielak and Orlicz
\cite{M,MO} ($f\in L^{\Phi}(\Omega)$ if and only if
$\int\Phi(x,|f(x)|)\;dx<\infty$). Many properties of Sobolev
spaces have been extended to Orlicz-Sobolev spaces, mainly by
Dankert \cite{dank}, Donaldson and Trudinger \cite{DT}, and
O'Neill \cite{onei} (see also Adams \cite{A} for an excellent
account of those works). Orlicz-Sobolev spaces have been used in
the last decades to model various phenomena. Chen, Levine and Rao
\cite{CLR} proposed a framework for image restoration based on a
variable exponent Laplacian. A second application which uses
variable exponent type Laplace operators is modelling
electrorheological fluids \cite{AM,R}.
 According to Diening
\cite{die}, we are strongly convinced that these more general
spaces will become increasingly  important in modelling modern
materials.
\smallskip

In this paper we assume that the function $a(x,t):\overline\Omega\times\RR\rightarrow\RR$ in \eq{1} is such that
$\phi(x,t):\overline\Omega\times\RR\rightarrow\RR$,
$$\phi(x,t)=\left\{\begin{array}{lll}
a(x,|t|)t, &\mbox{for}&
t\neq 0\\
0, &\mbox{for}& t=0\,,
\end{array}\right.$$
and satisfies
\smallskip

\noindent ($\phi$) for all $x\in\Omega$, $\phi(x,\cdot):\RR\rightarrow\RR$ is an odd, increasing homeomorphism from
$\RR$ onto $\RR$;
\smallskip

\noindent and $\Phi(x,t):\overline\Omega\times\RR\rightarrow\RR$,
$$\Phi(x,t)=\int_0^t\phi(x,s)\;ds,\;\;\;\forall\;x\in\overline\Omega,\;t\geq 0\,,$$
belongs to {\it class $\Phi$} (see \cite{M}, p. 33), i.e. $\Phi$ satisfies the following conditions
\smallskip

\noindent ($\Phi_1$) for all $x\in\Omega$, $\Phi(x,\cdot):[0,\infty)\rightarrow\RR$ is a nondecreasing continuous 
function,
with $\Phi(x,0)=0$ and $\Phi(x,t)>0$ whenever $t>0$; $\lim_{t\rightarrow\infty}\Phi(x,t)=\infty$;
\smallskip

\noindent ($\Phi_2$) for every $t\geq 0$, $\Phi(\cdot,t):\Omega\rightarrow\RR$ is a measurable function.
\medskip

\noindent{\bf Remark 1.} Since $\phi(x,\cdot)$ satisfies condition ($\phi$) we deduce that $\Phi(x,\cdot)$ is convex 
and
increasing from $\RR^+$ to $\RR^+$.

For the function $\Phi$ introduced above we define the {\it generalized Orlicz class},
$$K_\Phi(\Omega)=\{u:\Omega\rightarrow\RR,\;{\rm measurable};\;
\int_\Omega\Phi(x,|u(x)|)\;dx<\infty\}$$
and the {\it generalized Orlicz space},
$$L^\Phi(\Omega)=\{u:\Omega\rightarrow\RR,\;{\rm measurable};\;\lim_{\lambda\rightarrow 
0^+}\int_\Omega\Phi(x,\lambda|u(x)|)\;dx=0\}\,.$$
The space $L^\Phi(\Omega)$ is a Banach space endowed with the
{\it Luxemburg norm}
$$|u|_\Phi=\inf\left\{\mu>0;\ \int_\Omega\Phi\left(x,\frac{|u(x)|}{\mu}
\right)\;dx\leq 1\right\}$$
or the equivalent norm (the {\it Orlicz norm})
$$|u|_{(\Phi)}=\sup\left\{\left|\int_\Omega uv\;dx\right|;\ v\in L^{\overline\Phi}(\Omega),\
\int_\Omega \overline\Phi (x,|v(x)|)\;dx\leq 1\right\}\,,$$
where $\overline\Phi$ denotes the {\it conjugate Young} function of $\Phi$, that is,
$$\overline \Phi (x,t)=\sup_{s>0}\{ts-\Phi (x,s);\ s\in\RR\},\;\;\;\forall\;x\in\overline\Omega,\;t\geq 0\,.$$
Furthermore, for $\Phi$ and $\overline\Phi$ conjugate Young functions, the H\"older type inequality holds true
\begin{equation}\label{hol}
\left|\int_\Omega uv\;dx\right|\leq C\cdot|u|_\Phi\cdot|v|_{\overline\Phi},
\;\;\;\forall\;u\in L^\Phi(\Omega),\;v\in L^{\overline\Phi}(\Omega)\,,
\end{equation}
where $C$ is a positive constant (see \cite{M}, Theorem 13.13).

In this paper we assume that there exist two positive constants $\phi_0$ and $\phi^0$ such that
\begin{equation}\label{5}
1<\phi_0\leq\frac{t\phi(x,t)}{\Phi(x,t)}\leq\phi^0<\infty,\;\;\;\forall\;x\in\overline\Omega,\;t\geq 0\,.
\end{equation}
The above relation implies that $\Phi$ satisfies the
{\it$\Delta_2$-condition} (see Proposition \ref{p3}), i.e.
\begin{equation}\label{6}
\Phi(x,2t)\leq K\cdot\Phi(x,t),\;\;\;\forall\;x\in\overline\Omega,\;t\geq 0\,,
\end{equation}
where $K$ is a positive constant. Relation \eq{6} and Theorem 8.13 in \cite{M} imply that 
$L^\Phi(\Omega)=K_\Phi(\Omega)$.

Furthermore, we assume that $\Phi$ satisfies the following condition
\begin{equation}\label{7}
{\rm for}\;{\rm each}\;x\in\overline\Omega,\;{\rm the}\;{\rm function}\;[0,\infty)\ni 
t\rightarrow\Phi(x,\sqrt{t})\;{\rm is}\;{\rm convex}\,.
\end{equation}
Relation \eq{7} assures that $L^\Phi(\Omega)$ is an uniformly
convex space and thus, a reflexive space (see Proposition
\ref{p2}).

On the other hand, we point out that assuming that $\Phi$ and $\Psi$ belong to class $\Phi$ and
\begin{equation}\label{8}
\Psi(x,t)\leq K_1\cdot\Phi(x,K_2\cdot t)+h(x),\;\;\;\forall\;x\in\overline\Omega,\;t\geq 0\,,
\end{equation}
where $h\in L^1(\Omega)$, $h(x)\geq 0$ a.e. $x\in\Omega$ and $K_1$, $K_2$ are positive constants, then by Theorem 8.5 
in \cite{M}
we have that there exists the continuous embedding $L^\Phi(\Omega)\subset L^\Psi(\Omega)$.

An important role in manipulating the generalized Lebesgue-Sobolev
spaces
is played by the {\it modular} of the $L^{\Phi}(\Omega)$ space, which
is
the mapping
 $\rho_{\Phi}:L^{\Phi}(\Omega)\rightarrow\RR$ defined by
$$\rho_{\Phi}(u)=\int_\Omega\Phi(x,|u(x)|)\;dx.$$
If $(u_n)$, $u\in L^{\Phi}(\Omega)$  then the following relations hold true
\begin{equation}\label{9}
|u|_{\Phi}>1\;\;\;\Rightarrow\;\;\;|u|_{\Phi}^{\phi_0}\leq\rho_{\Phi}(u)
\leq|u|_{\Phi}^{\phi^0}\,,
\end{equation}
\begin{equation}\label{10}
|u|_{\Phi}<1\;\;\;\Rightarrow\;\;\;|u|_{\Phi}^{\phi^0}\leq
\rho_{\Phi}(u)\leq|u|_{\Phi}^{\phi_0}\,,
\end{equation}
\begin{equation}\label{11}
|u_n-u|_{\Phi}\rightarrow 0\;\;\;\Leftrightarrow\;\;\;\rho_{\Phi}
(u_n-u)\rightarrow 0\,,
\end{equation}
\begin{equation}\label{12}
|u_n|_{\Phi}\rightarrow\infty\;\;\;\Leftrightarrow\;\;\;\rho_{\Phi}
(u_n)\rightarrow\infty\,.
\end{equation}
Next, we define the {\it generalized Orlicz-Sobolev space}
$$W^{1,\Phi}(\Omega)=\left\{u\in L^\Phi(\Omega);\;\frac{\partial u}
{\partial x_i}\in L^\Phi(\Omega),\;i=1,...,N\right\}.$$
On $W^{1,\Phi}(\Omega)$ we define the equivalent norms
\begin{eqnarray*}
\|u\|_{1,\Phi}&=&|\;|\nabla u|\;|_\Phi+|u|_\Phi\\
\|u\|_{2,\Phi}&=&\max\{|\;|\nabla u|\;|_\Phi,|u|_\Phi\}\\
\|u\|&=&\inf\left\{\mu>0;\ \int_\Omega\left[\Phi\left(x,\frac{|u(x)|}{\mu}\right)+\Phi\left(x,\frac{|\nabla 
u(x)|}{\mu}
\right)\right]\;dx\leq 1\right\}\,,
\end{eqnarray*}
(see Proposition \ref{p4}).

The generalized Orlicz-Sobolev space $W^{1,\Phi}(\Omega)$ endowed with one of the above norms is a reflexive Banach 
space.

Finally, we point out that assuming that $\Phi$ and $\Psi$ belong to class $\Phi$, satisfying relation \eq{8} and
$\inf_{x\in\Omega}\Phi(x,1)>0$, $\inf_{x\in\Omega}\Psi(x,1)>0$ then there exists the continuous embedding
$W^{1,\Phi}(\Omega)\subset W^{1,\Psi}(\Omega)$.

We refer to Orlicz \cite{orl}, Nakano \cite{nak},  Musielak \cite{M}, Musielak and Orlicz \cite{MO}, Diening 
\cite{die} for further
properties of generalized Lebesgue-Sobolev spaces.
\medskip

\noindent{\bf Remark 2.} a) Assuming $\Phi(x,t)=\Phi(t)$, i.e. $\Phi$ is independent of variable $x$, we say that
$L^\Phi$ and $W^{1,\Phi}$ are {\it Orlicz spaces}, respectively {\it Orlicz-Sobolev spaces} (see 
\cite{A,Clem1,Clem2,orl}).
\smallskip

\noindent b) Assuming $\Phi(x,t)=|t|^{p(x)}$ with $p(x)\in C(\overline\Omega)$, $p(x)>1$ for all 
$x\in\overline\Omega$ we
denote  $L^\Phi$ by $L^{p(x)}$ and $W^{1,\Phi}$ by $W^{1,p(x)}$ and we refer to them as {\it variable exponents 
Lebesgue
spaces}, respectively {\it variable exponents Sobolev spaces} (see 
\cite{edm,edm2,edm3,FSZ,FZ1,KR,RoyalSoc,mihradproc,M,MO,nak})
\smallskip

\noindent c) Our framework enables us to work with spaces which are more general than those described in a) and b) 
(see the examples at the
end of this paper).

\section{Auxiliary results regarding generalized Orlicz-Sobolev spaces}
In this section we point out certain useful results regarding the generalized Orlicz-Sobolev spaces.

\begin{prop}\label{p1}
Assume condition \eq{5} is satisfied. Then the following relations hold true
\begin{equation}\label{13}
|u|_{\Phi}^{\phi_0}\leq\rho_{\Phi}(u)\leq|u|_{\Phi}^{\phi^0},\;\;\;\forall\;u\in L^\Phi(\Omega)\;{\rm 
with}\;|u|_\Phi>1\,,
\end{equation}
\begin{equation}\label{14}
|u|_{\Phi}^{\phi^0}\leq\rho_{\Phi}(u)\leq|u|_{\Phi}^{\phi_0},\;\;\;\forall\;u\in L^\Phi(\Omega)\;{\rm 
with}\;|u|_\Phi<1\,.
\end{equation}
\end{prop}
\proof First, we show that $\rho_{\Phi}(u)\leq|u|_{\Phi}^{\phi^0}$ for all $u\in L^\Phi(\Omega)$  with $|u|_\Phi>1$.

Indeed, since $\phi^0\geq(t\phi(x,t))/\Phi(x,t)$ for all $x\in\overline\Omega$ and all $t\geq 0$ it follows that
letting $\sigma>1$ we have
$$\log(\Phi(x,\sigma\cdot t))-\log(\Phi(x,t))=\int_t^{\sigma\cdot t}\frac{\phi(x,s)}{\Phi(x,s)}\;ds\leq
\int_t^{\sigma\cdot t}\frac{\phi^0}{s}\;ds=\log(\sigma^{\phi^0})\,.$$
Thus, we deduce
\begin{equation}\label{15}
\Phi(x,\sigma\cdot t)\leq\sigma^{\phi^0}\cdot\Phi(x,t),\;\;\;\forall\;x\in\overline\Omega,\;t>0,\;\sigma>1\,.
\end{equation}
Let now $u\in L^\Phi(\Omega)$ with $|u|_\Phi>1$. Using the definition of the Luxemburg norm and relation \eq{15} we 
deduce
\begin{eqnarray*}
\int_\Omega\Phi(x,|u(x)|)\;dx&=&\int_\Omega \Phi\left(x,|u|_\Phi\cdot\frac{|u(x)|}{|u|_\Phi}\right)\;dx\\
&\leq&|u|_{\Phi}^{\phi^0}\cdot\int_\Omega \Phi\left(x,\frac{|u(x)|}{|u|_\Phi}\right)\;dx\\
&\leq&|u|_{\Phi}^{\phi^0}\,.
\end{eqnarray*}
Now, we show that $\rho_{\Phi}(u)\geq|u|_{\Phi}^{\phi_0}$ for all $u\in L^\Phi(\Omega)$  with $|u|_\Phi>1$.

Since $\phi_0\leq(t\phi(x,t))/\Phi(x,t)$ for all $x\in\overline\Omega$ and all $t\geq 0$, similar techniques as those 
used
in the proof of relation \eq{15} imply
\begin{equation}\label{16}
\Phi(x,\sigma\cdot t)\geq\sigma^{\phi_0}\cdot\Phi(x,t),\;\;\;\forall\;x\in\overline\Omega,\;t>0,\;\sigma>1\,.
\end{equation}
Let  $u\in L^\Phi(\Omega)$ with $|u|_\Phi>1$. We consider $\beta\in(1,|u|_\Phi)$. Since $\beta<|u|_\Phi$ it follows 
that
$\int_\Omega \Phi\left(x,\frac{|u(x)|}{\beta}\right)\;dx>1$ otherwise we will obtain a contradiction with the 
definition
of the Luxemburg norm. The above considerations implies
$$\int_\Omega \Phi(x,|u(x)|)\;dx=\int_\Omega \Phi\left(x,\beta\cdot\frac{|u(x)|}{\beta}\right)\;dx\geq
\beta^{\phi_0}\cdot\int_\Omega \Phi\left(x,\frac{|u(x)|}{\beta}\right)\;dx\geq\beta^{\phi_0}\,.$$
Letting $\beta\nearrow|u|_\Phi$ we deduce that relation \eq{13} holds true.

Next, we show that $\rho_{\Phi}(u)\leq|u|_{\Phi}^{\phi_0}$ for all $u\in L^\Phi(\Omega)$  with $|u|_\Phi<1$.
It is easy to show (see the proof of relations \eq{15} and \eq{16}) that
\begin{equation}\label{17}
\Phi(x,t)\leq\tau^{\phi_0}\cdot\Phi(x,t/\tau),\;\;\;\forall\;x\in\overline\Omega,\;t>0,\;\tau\in(0,1)\,.
\end{equation}
Let  $u\in L^\Phi(\Omega)$ with $|u|_\Phi<1$. The definition of the Luxemburg norm and relation \eq{17} imply
\begin{eqnarray*}
\int_\Omega\Phi(x,|u(x)|)\;dx&=&\int_\Omega \Phi\left(x,|u|_\Phi\cdot\frac{|u(x)|}{|u|_\Phi}\right)\;dx\\
&\leq&|u|_{\Phi}^{\phi_0}\cdot\int_\Omega \Phi\left(x,\frac{|u(x)|}{|u|_\Phi}\right)\;dx\\
&\leq&|u|_{\Phi}^{\phi_0}\,.
\end{eqnarray*}

Finally, we show that $\rho_{\Phi}(u)\geq|u|_{\Phi}^{\phi^0}$ for all $u\in L^\Phi(\Omega)$  with $|u|_\Phi<1$.

As in the proof of \eq{15} we deduce
\begin{equation}\label{18}
\Phi(x,t)\geq\tau^{\phi^0}\cdot\Phi(x,t/\tau),\;\;\;\forall\;x\in\overline\Omega,\;t>0,\;\tau\in(0,1)\,.
\end{equation}
Let  $u\in L^\Phi(\Omega)$ with $|u|_\Phi<1$ and $\xi\in(0,|u|_\Phi)$. By \eq{18} we find
\begin{equation}\label{19}
\int_\Omega \Phi(x,|u(x)|)\;dx\geq\xi^{\phi^0}\cdot\int_\Omega \Phi\left(x,\frac{|u(x)|}{\xi}\right)\;dx\,.
\end{equation}
Define $v(x)=u(x)/\xi$, for all $x\in\Omega$. We have $|v|_\Phi=|u|_\Phi/\xi>1$. Using relation \eq{13} we find
\begin{equation}\label{20}
\int_\Omega \Phi(x,|v(x)|)\;dx\geq|v|_\Phi^{\phi_0}>1\,.
\end{equation}
By \eq{19} and \eq{20} we obtain
$$\int_\Omega \Phi(x,|u(x)|)\;dx\geq\xi^{\phi^0},\;\;\;\forall\;\xi\in(0,|u|_\Phi)\,.$$
Letting $\xi\nearrow|u|_\Phi$ we deduce that relation \eq{14} holds true. The proof of Proposition \ref{p1} is 
complete.  \endproof

\begin{prop}\label{p2}
Assume $\Phi$ satisfies conditions \eq{6} and \eq{7}. Then the space $L^\Phi(\Omega)$ is uniformly convex.
\end{prop}
\proof
From the above hypotheses we deduce that we can apply Lemma 2.1 in \cite{Lam} in order to deduce
$$\frac{1}{2}[\Phi(x,|t|)+\Phi(x,|s|)]\geq\Phi\left(x,\frac{|t+s|}{2}\right)+\Phi\left(x,\frac{|t-s|}{2}\right),\;\;\%
;
\forall\;x\in\Omega,\;t,\;s\in\RR\,.$$
The above inequality yields
\begin{equation}\label{L}
\frac{1}{2}[\rho_\Phi(u)+\rho_\Phi(v)]\geq\rho_\Phi\left(\frac{u+v}{2}\right)+\rho_\Phi\left(\frac{u-v}{2}\right),\;\%
;\;
\forall\;u,v\in L^\Phi(\Omega)\,.
\end{equation}
Assume that $|u|_\Phi<1$ and $|v|_\Phi<1$ and $|u-v|_\Phi>\epsilon$ (with $\epsilon\in(0,1/K)$). Then we have
\begin{eqnarray*}
\rho_\Phi(u-v)&\geq&|u-v|_\Phi^{\phi_0}\;\;\;{\rm if}\;|u-v|_\Phi>1\\
\rho_\Phi(u-v)&\geq&|u-v|_\Phi^{\phi^0}\;\;\;{\rm if}\;|u-v|_\Phi<1\,,
\end{eqnarray*}
and
$$\rho_\Phi(u)<1,\;\;\; \rho_\Phi(v)<1\,.$$
The above information and relation \eq{5} yield
$$\rho_\Phi\left(\frac{u-v}{2}\right)\geq\frac{1}{K}\cdot\rho_\Phi(u-v)\geq
\left\{\begin{array}{lll}
\di\frac{1}{K}\cdot\epsilon^{\phi_0}, &\mbox{if}&
|u-v|_\Phi>1\\
\di\frac{1}{K}\cdot\epsilon^{\phi^0}, &\mbox{if}&
|u-v|_\Phi<1\,.
\end{array}\right.$$
By \eq{L} and the above inequality we have
\begin{equation}\label{L1}
\rho_\Phi\left(\frac{u+v}{2}\right)<
\left\{\begin{array}{lll}
1-\di\frac{1}{K}\cdot\epsilon^{\phi_0}, &\mbox{if}&
|u-v|_\Phi>1\\
1-\di\frac{1}{K}\cdot\epsilon^{\phi^0}, &\mbox{if}&
|u-v|_\Phi<1\,.
\end{array}\right.
\end{equation}
On the other hand, we have
\begin{equation}\label{L2}
\rho_\Phi\left(\frac{u+v}{2}\right)\geq
\left\{\begin{array}{lll}
\di\left|\frac{u+v}{2}\right|_\Phi^{\phi_0}, &\mbox{if}&
\di\left|\frac{u+v}{2}\right|_\Phi>1\\
\di\left|\frac{u+v}{2}\right|_\Phi^{\phi^0}, &\mbox{if}&
\di\left|\frac{u+v}{2}\right|_\Phi<1\,.
\end{array}\right.
\end{equation}
Relations \eq{L1} and \eq{L2} show that there exists $\delta>0$ such that
$$\left|\frac{u+v}{2}\right|_\Phi<1-\delta\,.$$
Thus, we proved that $L^\Phi(\Omega)$ is an uniformly convex space. The proof of Proposition \ref{p2} is complete.  
\endproof

\noindent{\bf Remark 3.} Condition \eq{7} (via relation \eq{L}) also implies the fact that for every $x\in\Omega$ 
fixed, the
function  $\Phi(x,\cdot)$ is convex from $\RR^+$ to $\RR^+$.

\begin{prop}\label{p3}
Condition \eq{5} implies condition \eq{6}.
\end{prop}
\proof
Since relation \eq{5} holds true by Proposition \ref{p1} it follows that condition \eq{15} works. We deduce that
$$\Phi(x,2\cdot t)\leq 2^{\phi^0}\cdot\Phi(x,t),\;\;\;\forall\;x\in\Omega,\;t>0\,.$$
Thus, relation \eq{6} holds true with $K=2^{\phi^0}$. The proof of Proposition \ref{p3} is complete.  \endproof

\begin{prop}\label{p4}
On $W^{1,\Phi}(\Omega)$ the following norms
$$\|u\|_{1,\Phi}=|\;|\nabla u|\;|_\Phi+|u|_\Phi\,,$$
$$\|u\|_{2,\Phi}=\max\{|\;|\nabla u|\;|_\Phi,|u|_\Phi\}\,,$$
$$\|u\|=\inf\left\{\mu>0;\ \int_\Omega\left[\Phi\left(x,\frac{|u(x)|}{\mu}
\right)+\Phi\left(x,\frac{|\nabla u(x)|}{\mu}\right)\right]\;dx\leq 1\right\}\,,$$
are equivalent.
\end{prop}
\proof First, we point out that  $\|\;\|_{1,\Phi}$ and $\|\;\|_{2,\Phi}$ are equivalent, since
\begin{equation}\label{20biss}
2\cdot\|u\|_{2,\Phi}\geq\|u\|_{1,\Phi}\geq\|u\|_{2,\Phi},\;\;\;\forall\;u\in W^{1,\Phi}(\Omega)\,.
\end{equation}
Next, we remark that
$$
\int_\Omega\Phi\left(x,\frac{|u(x)|}{|u|_\Phi}
\right)\;dx\leq 1\;\;\;{\rm and}\;\;\;\int_\Omega\Phi\left(x,\frac{|\nabla u(x)|}{|\;|\nabla u|\;|_\Phi}
\right)\;dx\leq 1\,,
$$
and
$$
\int_\Omega\left[\Phi\left(x,\frac{|u(x)|}{\|u\|}
\right)+\Phi\left(x,\frac{|\nabla u(x)|}{\|u\|}\right)\right]\;dx\leq 1\,.
$$
Using the above relations we obtain
$$\int_\Omega\Phi\left(x,\frac{|u(x)|}{\|u\|}
\right)\;dx\leq 1\;\;\;{\rm and}\;\;\;\int_\Omega\Phi\left(x,\frac{|\nabla u(x)|}{\|u\|}
\right)\;dx\leq 1.$$
Taking into account the way in which $|\;|_\Phi$ is defined we find
\begin{equation}\label{21}
2\|u\|\geq(|u|_\Phi+|\;|\nabla u|\;|_\Phi)=\|u\|_{1,\Phi},\;\;\;\forall\;u\in W^{1,\Phi}(\Omega)\,.
\end{equation}
On the other hand, by relation \eq{16} we deduce that
$$\Phi(x,2\cdot t)\geq 2\cdot\Phi(x,t),\;\;\;\forall\;x\in\Omega,\;t>0\,.$$
Thus, we deduce that
$$2\cdot\Phi\left(x,\frac{|u(x)|}{2\cdot\|u\|_{2,\Phi}}\right)\leq\Phi\left(x,\frac{|u(x)|}{\|u\|_{2,\Phi}}\right),\;
\;\;\forall\;
u\in W^{1,\Phi}(\Omega),\;x\in\Omega$$
 and
$$2\cdot\Phi\left(x,\frac{|\nabla u(x)|}{2\cdot\|u\|_{2,\Phi}}\right)\leq\Phi\left(x,\frac{|\nabla 
u(x)|}{\|u\|_{2,\Phi}}\right),\;\;\;\forall\;
u\in W^{1,\Phi}(\Omega),\;x\in\Omega\,.$$
It follows that
\begin{equation}\label{22}
\int_\Omega\left[\Phi\left(x,\frac{|u(x)|}{2\|u\|_{2,\Phi}}
\right)+\Phi\left(x,\frac{|\nabla u(x)|}{2\|u\|_{2,\Phi}}\right)\right]\;dx\leq\frac{1}{2}
\left\{\int_\Omega\left[\Phi\left(x,\frac{|u(x)|}{\|u\|_{2,\Phi}}
\right)+\Phi\left(x,\frac{|\nabla u(x)|}{\|u\|_{2,\Phi}}\right)\right]\;dx\right\}.
\end{equation}
But, since
$$\|u\|_{2,\Phi}\geq|u|_\Phi\;\;\;{\rm and}\;\;\;\|u\|_{2,\Phi}\geq|\;|\nabla u|\;|_\Phi,\;\;\;\forall\;u\in 
W^{1,\Phi}(\Omega)\,,$$
we obtain
\begin{equation}\label{23}
\frac{|u(x)|}{|u|_\Phi}\geq\frac{|u(x)|}{\|u\|_{2,\Phi}}\;\;\;{\rm and}\;\;\;
\frac{|\nabla u(x)|}{|\;|\nabla u|\;|_\Phi}\geq\frac{|\nabla u(x)|}{\|u\|_{2,\Phi}},\;\;\;\forall\;u\in 
W^{1,\Phi}(\Omega),\;x\in\Omega\,.
\end{equation}
Taking into account that $\Phi$ is increasing by \eq{22} and
\eq{23} we deduce that
$$\int_\Omega\left[\Phi\left(x,\frac{|u(x)|}{2\|u\|_{2,\Phi}}
\right)+\Phi\left(x,\frac{|\nabla u(x)|}{2\|u\|_{2,\Phi}}\right)\right]\;dx\leq\frac{1}{2}
\left\{\int_\Omega\left[\Phi\left(x,\frac{|u(x)|}{|u|_{\Phi}}
\right)+\Phi\left(x,\frac{|\nabla u(x)|}{|\;|\nabla u\;|_{\Phi}}\right)\right]\;dx\right\}\leq 1\,,$$
for all $u\in W^{1,\Phi}(\Omega)$.

We conclude that
\begin{equation}\label{24}
2\cdot\|u\|_{1,\Phi}\geq 2\cdot\|u\|_{2,\Phi}\geq\|u\|,\;\;\;\forall\;u\in W^{1,\Phi}(\Omega)\,.
\end{equation}
By relations \eq{20biss}, \eq{21} and \eq{24} we deduce that Proposition \ref{p4} holds true.  \endproof

\begin{prop}\label{p5}
The following relations hold true
\begin{equation}\label{25}
\int_\Omega[\Phi(x,|u(x)|)+\Phi(x,|\nabla u(x)|)]\;dx\geq\|u\|^{\phi_0},\;\;\;\forall\;u\in W^{1,\Phi}(\Omega)\;{\rm 
with}\;\|u\|>1\,;
\end{equation}
\begin{equation}\label{26}
\int_\Omega[\Phi(x,|u(x)|)+\Phi(x,|\nabla u(x)|)]\;dx\geq\|u\|^{\phi^0},\;\;\;\forall\;u\in W^{1,\Phi}(\Omega)\;{\rm 
with}\;\|u\|<1\,.
\end{equation}
\end{prop}
\proof
First, assume that $\|u\|>1$. Let $\beta\in(1,\|u\|)$. By relation \eq{16} we have
$$\int_\Omega[\Phi(x,|u(x)|)+\Phi(x,|\nabla 
u(x)|)]\;dx\geq\beta^{\phi_0}\cdot\int_\Omega\left[\Phi\left(x,\frac{|u(x)|}{\beta}\right)+
\Phi\left(x,\frac{|\nabla u(x)|}{\beta}\right)\right]\,.$$
Since $\beta<\|u\|$ we find
$$\int_\Omega\left[\Phi\left(x,\frac{|u(x)|}{\beta}\right)+
\Phi\left(x,\frac{|\nabla u(x)|}{\beta}\right)\right]>1\,.$$
Thus, we find
$$\int_\Omega[\Phi(x,|u(x)|)+\Phi(x,|\nabla u(x)|)]\;dx\geq\beta^{\phi_0}\,.$$
Letting $\beta\nearrow\|u\|$ we deduce that \eq{25} holds true.

Next, assume $\|u\|<1$. Let $\xi\in(0,\|u\|)$. By relation \eq{18} we obtain
\begin{equation}\label{27}
\int_\Omega[\Phi(x,|u(x)|)+\Phi(x,|\nabla 
u(x)|)]\;dx\geq\xi^{\phi^0}\cdot\int_\Omega\left[\Phi\left(x,\frac{|u(x)|}{\xi}\right)+
\Phi\left(x,\frac{|\nabla u(x)|}{\xi}\right)\right]\;dx.
\end{equation}
Defining $v(x)=u(x)/\xi$, for all $x\in\Omega$, we have
$\|v\|=\|u\|/\xi>1$. Using relation \eq{25} we find
\begin{equation}\label{28}
\int_\Omega[\Phi(x,|v(x)|)+\Phi(x,|\nabla v(x)|)]\;dx\geq\|v\|^{\phi_0}>1.
\end{equation}
Relations \eq{27} and \eq{28} show that
$$\int_\Omega[\Phi(|u(x)|)+\Phi(|\nabla u(x)|)]\;dx\geq\xi^{p^0}.$$
Letting $\xi\nearrow\|u\|$ in the above inequality we obtain that relation \eq{26} holds true. The proof of
Proposition \ref{p5} is complete.  \endproof

\section{Main results}
In this paper we study problem \eq{1} in the particular case when $\Phi$ satisfies
\begin{equation}\label{29}
M\cdot|t|^{p(x)}\leq\Phi(x,t),\;\;\;\forall\;x\in\overline\Omega,\;t\geq 0\,,
\end{equation}
where $p(x)\in C(\overline\Omega)$ with $p(x)>1$ for all $x\in\overline\Omega$ and $M>0$ is a constant.
\smallskip

\noindent{\bf Remark 4.} By relation \eq{29} we deduce that
$W^{1,\Phi}(\Omega)$ is continuously embedded in
$W^{1,p(x)}(\Omega)$ (see relation \eq{8} with
$\Psi(x,t)=|t|^{p(x)}$). On the other hand, it is known (see
\cite{KR,FSZ,RoyalSoc}) that  $W^{1,p(x)}(\Omega)$ is compactly
embedded in $L^{r(x)}(\Omega)$ for any $r(x)\in
C(\overline\Omega)$ with $1<r^-\leq r^+<\frac{Np^-}{N-p^-}$. Thus,
we deduce that $W^{1,\Phi}(\Omega)$ is compactly embedded in
$L^{r(x)}(\Omega)$ for any $r(x)\in C(\overline\Omega)$ with
$1<r(x)<\frac{Np^-}{N-p^-}$ for all $x\in\overline\Omega$.

On the other hand, we assume that the function $g$ from problem \eq{1} satisfies the hypotheses
\begin{equation}\label{g}
|g(x,t)|\leq C_0\cdot|t|^{q(x)-1},\;\;\;\forall\;x\in\Omega,\;t\in\RR
\end{equation}
and
\begin{equation}\label{G1}
C_1\cdot|t|^{q(x)}\leq G(x,t):=\int_0^tg(x,s)\;ds\leq C_2\cdot|t|^{q(x)}, \;\;\;\forall\;x\in\Omega,\;t\in\RR\,,
\end{equation}
where $C_0$, $C_1$ and $C_2$ are positive constants and $q(x)\in C(\overline\Omega)$ satisfies
$1<q(x)<\frac{Np^-}{N-p^-}$ for all $x\in\overline\Omega$.
\smallskip

\noindent{\bf Examples.} We point out certain examples of functions $g$ and $G$ which satisfy hypotheses \eq{g} and 
\eq{G1}.

1) $g(x,t)=q(x)\cdot|t|^{q(x)-2}t$ and $G(x,t)=|t|^{q(x)}$, where $q(x)\in C(\overline\Omega)$ satisfies
$2\leq q(x)<\frac{Np^-}{N-p^-}$ for all $x\in\overline\Omega$;

2) $g(x,t)=q(x)\cdot|t|^{q(x)-2}t+(q(x)-2)\cdot[\log(1+t^2)]\cdot|t|^{q(x)-4}t+\frac{t}{1+t^2}|t|^{q(x)-2}$ and
$G(x,t)=|t|^{q(x)}+\log(1+t^2)\cdot|t|^{q(x)-2}$, where $q(x)\in C(\overline\Omega)$ satisfies
$4\leq q(x)<\frac{Np^-}{N-p^-}$ for all $x\in\overline\Omega$;

3) $g(x,t)=q(x)\cdot|t|^{q(x)-2}t+(q(x)-1)\cdot 
\sin(\sin(t))\cdot|t|^{q(x)-3}t+\cos(\sin(t))\cdot\cos(t)\cdot|t|^{q(x)-1}$ and
$G(x,t)=|t|^{q(x)}+\sin(\sin(t))\cdot|t|^{q(x)-1}$, where $q(x)\in C(\overline\Omega)$ satisfies
$3\leq q(x)<\frac{Np^-}{N-p^-}$ for all $x\in\overline\Omega$.
\smallskip

We say that $u\in W^{1,\Phi}(\Omega)$ is a {\it weak solution} of problem \eq{1} if
$$\int_{\Omega}a(x,|\nabla u|)\nabla u\nabla v\;dx+\int_\Omega a(x,|u|)uv\;dx-\lambda\int_{\Omega}g(x,u)v\;dx=0,$$
for all $v\in W^{1,\Phi}(\Omega)$.

The main results of this paper are given by the following theorems.
\begin{teo}\label{t1}
Assume $\phi$ and $\Phi$ verify conditions ($\phi$), ($\Phi_1$), ($\Phi_2$), \eq{5}, \eq{7} and \eq{29} and the 
functions $g$ and
$G$ satisfy conditions \eq{g} and \eq{G1}. Furthermore, we assume that $q^-<\phi_0$. Then there exists 
$\lambda_\star>0$ such that
for any $\lambda\in(0,\lambda_\star)$ problem \eq{1} has a nontrivial weak solution.
\end{teo}
\begin{teo}\label{t2}
Assume $\phi$ and $\Phi$ verify conditions ($\phi$), ($\Phi_1$), ($\Phi_2$), \eq{5}, \eq{7} and \eq{29} and the 
functions $g$ and
$G$ satisfy conditions \eq{g} and \eq{G1}. Furthermore, we assume that $q^+<\phi_0$. Then there exists 
$\lambda_\star>0$ and
$\lambda^\star>0$ such that
for any $\lambda\in(0,\lambda_\star)\cup(\lambda^\star,\infty)$ problem \eq{1} has a nontrivial weak solution.
\end{teo}

\section{Proof of the main results}
Let $E$ denote the generalized Orlicz-Sobolev space $W^{1,\Phi}(\Omega)$.

For each $\lambda>0$ we define the energy functional $J_\lambda:E\rightarrow\RR$ by
$$J_\lambda(u)=\int_\Omega[\Phi(x,|\nabla u|)+\Phi(x,|u|)]\;dx-\lambda\int_\Omega G(x,u)\;dx,\;\;\;\forall\;u\in 
E\,.$$
We first establish some basic properties of $J_\lambda$.

\begin{prop}\label{p6}
For each $\lambda>0$ the functional $J_\lambda$ is well-defined on $E$ and $J_\lambda\in C^1(E,\RR)$ with
the derivative given by
$$\langle J_\lambda^{'}(u),v\rangle=\int_{\Omega}a(x,|\nabla u|)\nabla u\cdot\nabla v\;dx+\int_\Omega 
a(x,|u|)uv\;dx-\lambda
\int_\Omega g(x,u)v\;dx\,,$$
for all $u$, $v\in E$.
\end{prop}
To prove Proposition \ref{p6} we define the functional $\Lambda:E\rightarrow\RR$ by
$$\Lambda(u)=\int_\Omega[\Phi(x,|\nabla u|)+\Phi(x,|u|)]\;dx,\;\;\;\forall\;u\in E\,.$$
\begin{lemma}\label{l1}
The functional $\Lambda$ is well defined on $E$ and $\Lambda\in C^1(E,\RR)$ with
$$\langle\Lambda^{'}(u),v\rangle=\int_{\Omega}a(x,|\nabla u|)\nabla u\cdot\nabla v\;dx+\int_\Omega 
a(x,|u|)uv\;dx\,,$$
for all $u$, $v\in E$.
\end{lemma}
\proof
Clearly, $\Lambda$ is well defined on $E$.

{\bf Existence of the G\^ateaux derivative.} Let
$u$, $v\in E$. Fix $x\in\Omega$ and $0<|r|<1$. Then, by the mean value theorem, there exists $\nu,\theta\in[0,1]$ 
such that
\begin{equation}\label{30}
\begin{array}{lll}
|\Phi(x,|\nabla u(x)+r\nabla v(x)|)-\Phi(x,|\nabla u(x)|)|/|r|&=&
|\phi(x,|(1-\nu)|\nabla u(x)+r\nabla v(x)|+\nu|\nabla u(x)|)|\cdot\\
&&||\nabla u(x)+r\nabla v(x)|-|\nabla u(x)||\,
\end{array}
\end{equation}
and
\begin{equation}\label{31}
\begin{array}{lll}
|\Phi(x,| u(x)+r v(x)|)-\Phi(x,|u(x)|)|/|r|&=&
|\phi(x,|(1-\theta)| u(x)+r v(x)|+\theta| u(x)|)|\cdot\\
&&|| u(x)+r v(x)|-| u(x)||\,.
\end{array}
\end{equation}
Next, we claim that $\phi(x,|u(x)|)\in L^{\overline\Phi}(\Omega)$ provided that $u\in L^{\Phi}(\Omega)$, where 
$\overline\Phi$ is
the  conjugate Young function of $\Phi$.

Indeed, we know that
$$\overline \Phi (x,t)=\sup_{s>0}\{ts-\Phi (x,s);\ s\in\RR\},\;\;\;\forall\;x\in\overline\Omega,\;t\geq 0\,$$
or
$$\overline \Phi (x,t)=\int_0^t\overline\phi(x,s)\;ds,\;\;\;\forall\;x\in\overline\Omega,\;t\geq 0\,,$$
where $\overline\phi(x,t)=\sup\limits_{\phi(x,s)\leq t}s$, for all $x\in\overline\Omega$ and $t\geq 0$.

On the other hand, by relation ($\phi$) we know that for all $x\in\Omega$, $\phi(x,\cdot):\RR\rightarrow\RR$ is an
odd, increasing homeomorphism from $\RR$ onto $\RR$ and thus, an increasing homeomorphism from $\RR^+$ onto $\RR^+$.
It follows that for each $x\in\overline\Omega$ we can denote by $\phi^{-1}(x,t)$ the inverse function of $\phi(x,t)$
relative to  variable $t$. Thus, we deduce that $\phi(x,s)\leq t$ if and only if $s\leq\phi^{-1}(x,t)$. Taking into
account the above pieces of information we deduce that $\overline\phi(x,t)=\phi^{-1}(x,t)$. Consequently we find
$$\overline \Phi (x,t)=\int_0^t\phi^{-1}(x,s)\;ds,\;\;\;\forall\;x\in\overline\Omega,\;t\geq 0\,.$$
Next, since
$$\overline \Phi 
(x,\phi^{-1}(x,s))=\int_0^{\phi^{-1}(x,s)}\phi(x,\theta)\;d\theta,\;\;\;\forall\;x\in\overline\Omega,\;s\geq 0\,,$$
taking $\phi(x,\theta)=r$ we find
$$\overline \Phi (x,\phi^{-1}(x,s))=\int_0^sr\cdot(\phi^{-1}(x,r))^{'}_r\;dr=s\cdot\phi^{-1}(x,s)-\overline \Phi 
(x,s),
\;\;\;\forall\;x\in\overline\Omega,\;s\geq 0\,.$$
The above relation implies
$$\overline \Phi (x,s)\leq s\cdot\phi^{-1}(x,s),\;\;\;\forall\;x\in\overline\Omega,\;s\geq 0\,.$$
Taking into the above inequality $s=\phi(x,t)$ we find
$$\overline \Phi (x,\phi(x,t))\leq t\cdot\phi(x,t),\;\;\;\forall\;x\in\overline\Omega,\;t\geq 0\,.$$
The last inequality and relation \eq{5} yield
$$\overline \Phi (x,\phi(x,t))\leq \phi^0\cdot\Phi(x,t),\;\;\;\forall\;x\in\overline\Omega,\;t\geq 0\,.$$
Thus, we deduce that for any $u\in L^{\Phi}(\Omega)$ we have $\phi(x,|u(x)|)\in L^{\overline\Phi}(\Omega)$ and the 
claim is verified.
By applying relations \eq{30}, \eq{31}, the above claim and \eq{hol} we infer that
\begin{eqnarray*}
|\Phi(x,|\nabla u(x)+r\nabla v(x)|)+\Phi(x,| u(x)+r v(x)|)-\Phi(x,|\nabla u(x)|)-\Phi(x,|u(x)|)||/|r|\leq\\
|\phi(x,|(1-\nu)|\nabla u(x)+r\nabla v(x)|+\nu|\nabla u(x)|)|\cdot
||\nabla u(x)+r\nabla v(x)|-|\nabla u(x)||+\\|\phi(x,|(1-\theta)| u(x)+r v(x)|+\theta|\nabla u(x)|)|\cdot
|| u(x)+r v(x)|-| u(x)||\in L^1(\Omega)\,,
\end{eqnarray*}
for all $u$, $v\in E$, $x\in\overline\Omega$ and $|r|\in(0,1)$. It follows from the Lebesgue theorem that
$$\langle\Lambda^{'}(u),v\rangle=\int_{\Omega}a(x,|\nabla u|)\nabla u\cdot\nabla v\;dx+\int_\Omega 
a(x,|u|)uv\;dx\,.$$
{\bf Continuity of the G\^ateaux derivative.} Assume $u_n\rightarrow u$
in $E$. The above claim and the Lebesgue theorem imply
$$a(x,|\nabla u_n|)\nabla u_n\rightarrow a(x,|\nabla u|)\nabla u,\;\;\;{\rm in}\;\;\;(L^{\overline\Phi}(\Omega))^N$$
and
$$a(x,| u_n|) u_n\rightarrow a(x,| u|) u,\;\;\;{\rm in}\;\;\;L^{\overline\Phi}(\Omega)\,.$$
Those facts and \eq{hol} imply
$$|\langle\Lambda^{'}(u_n)-\Lambda^{'}(u),v\rangle|\leq
|a(x,|\nabla u_n|)\nabla u_n-a(x,|\nabla u|)\nabla u|_{\overline\Phi}\cdot|\;|\nabla v|\;|_\Phi+
|a(x,| u_n|) u_n-a(x,| u|) u|_{\overline\Phi}\cdot| v|_\Phi\,,$$
for all $v\in E$, and so
$$\|\Lambda^{'}(u_n)-\Lambda^{'}(u)\|\leq|a(x,|\nabla u_n|)\nabla u_n-a(x,|\nabla u|)\nabla u|_{\overline\Phi}+
|a(x,| u_n|) u_n-a(x,| u|) u|_{\overline\Phi}\rightarrow 0,\;\;\;{\rm as}\;{n\rightarrow\infty}\,.$$
The proof of Lemma \ref{l1} is complete.
\endproof

Combining Lemma \ref{l1} and Remark 4 we infer that Proposition \ref{p6} holds true.

\begin{lemma}\label{l2}
The functional $\Lambda$ is weakly lower semi-continuous.
\end{lemma}
\proof
By Corollary III.8 in \cite{B}, it is enough to show that $\Lambda$ is  lower semi-continuous. For this purpose, we 
fix
$u\in E$ and $\epsilon>0$. Since $\Lambda$ is convex (because $\Phi$ is convex) we deduce that for any $v\in E$ the 
following
inequality holds true
$$\Lambda(v)\geq\Lambda(u)+\langle\Lambda^{'}(u),v-u\rangle\,,$$
or
\begin{eqnarray*}
\Lambda(v)&\geq&\Lambda(u)-\int_\Omega[a(x,|\nabla u|)|\nabla u|\cdot|\nabla v-\nabla u|+
a(x,| u|)| u|\cdot| v- u|]\;dx\\
&=&\Lambda(u)-\int_\Omega[\phi(x,|\nabla u|)|\nabla v-\nabla u|+
\phi(x,| u|)| v- u|]\;dx\,.
\end{eqnarray*}
But, by the claim proved in Proposition \ref{p6} we know that for
any $u\in L^\Phi(\Omega)$ we have $\phi(x,|u|)$, $\phi(x,|\nabla
u|)\in L^{\overline\Phi}(\Omega)$. Thus, by relation \eq{hol} we
find
\begin{eqnarray*}
\Lambda(v)&\geq&\Lambda(u)-C\cdot[|\phi(x,|\nabla u|)|_{\overline\Phi}\cdot|\;|\nabla v-\nabla u|\;|_\Phi+
|\phi(x,| u|)|_{\overline\Phi}\cdot| v-u|_\Phi]\\
&\geq&\Lambda(u)-C^{'}\cdot\|u-v\|\\
&\geq&\Lambda(u)-\epsilon\,,
\end{eqnarray*}
for all $v\in E$ with $\|v-u\|<\delta=\epsilon/{C^{'}}$, where $C$ and
$C^{'}$ are positive constants. The proof
of Lemma \ref{l2} is complete.   \endproof

\begin{prop}\label{p7}
The functional $J_\lambda$ is weakly lower semi-continuous.
\end{prop}
\proof
Using Lemma \ref{l2} we have that $\Lambda$ is weakly lower semi-continuous. We show that $J_\lambda$ is weakly lower 
semi-continuous.
Let $\{u_n\}\subset E$ be a sequence which converges weakly to $u$ in $E$. By Lemma \ref{l2} we deduce
$$\Lambda(u)\leq\liminf_{n\rightarrow\infty}\Lambda(u_n)\,.$$
On the other hand, Remark 4 and conditions \eq{g} and \eq{G1} imply
$$\lim_{n\rightarrow\infty}\int_\Omega G(x,u_n)\;dx=\int_\Omega G(x,u)\;dx\,.$$
Thus, we find
$$J_\lambda(u)\leq\liminf_{n\rightarrow\infty}J_\lambda(u_n)\,.$$
Therefore, $J_\lambda$ is weakly lower semi-continuous and Proposition \ref{p7} is verified.   \endproof

\begin{prop}\label{p8}
Assume that the sequence $\{u_n\}$ converges weakly to $u$ in $E$ and
$$\limsup\limits_{n\rightarrow\infty}\langle\Lambda^{'}(u_n),u_n-u\rangle\leq 0.$$
Then $\{u_n\}$ converges strongly to $u$ in $E$.
\end{prop}
\proof
Since $\{u_n\}$ converges weakly to $u$ in $E$ it follows that $\{\|u_n\|\}$ is a bounded sequence of real numbers. 
That
fact and Proposition \ref{p4} imply that $\{|u_n|_\Phi\}$ and $\{|\;|\nabla u_n|\;|_\Phi\}$ are bounded sequences of 
real numbers.
That information and relations \eq{9} and \eq{10} yield that the sequence $\{\Lambda(u_n)\}$ is bounded. Then, up to 
a subsequence,
we deduce that $\Lambda(u_n)\rightarrow c$.

By Lemma \ref{l2} we obtain
$$\Lambda(u)\leq\liminf_{n\rightarrow\infty}\Lambda(u_n)=c\,.$$
On the other hand, since $\Lambda$ is convex, we have
$$\Lambda(u)\geq\Lambda(u_n)+\langle\Lambda^{'}(u_n),u-u_n\rangle\,.$$
Using the above hypothesis we conclude that $\Lambda(u)=c$. Taking into account that $\{(u_n+u)/2\}$ converges weakly 
to $u$ in
$E$ and using Lemma \ref{l2} we find
\begin{equation}\label{35}
c=\Lambda(u)\leq\Lambda\left(\frac{u_n+u}{2}\right)\,.
\end{equation}
We assume by contradiction that $\{u_n\}$ does not converge to $u$
in $E$ or $\{(u_n-u)/2\}$ does not converge to $0$ in $E$. It follows that there exist $\epsilon>0$ and
a subsequence $\{u_{n_m}\}$ of $\{u_n\}$ such that
\begin{equation}\label{36}
\left\|\frac{u_{n_m}-u}{2}\right\|\geq\epsilon,\;\;\;
\forall m\,.
\end{equation}
Furthermore, relations \eq{25}, \eq{26} and \eq{36} imply that there exists $\epsilon_1>0$ such that
\begin{equation}\label{37}
\Lambda\left(\frac{u_{n_m}-u}{2}\right)\geq\epsilon_1,\;\;\;
\forall m\,.
\end{equation}
On the other hand, relations \eq{L} and \eq{37} yield
$$\frac{1}{2}\Lambda(u)+\frac{1}{2}\Lambda(u_{n_m})-\Lambda\left(\frac{u_{n_m}+u}{2}\right)\geq\Lambda\left(\frac{u_{
n_m}-u}{2}\right)\geq
\epsilon_1,\;\;\;\forall m\,.$$
Letting $m\rightarrow\infty$ in the above inequality we obtain
$$c-\epsilon_1\geq\limsup_{m\rightarrow\infty}\Lambda\left(\frac{u_{n_m}+u}{2}\right)\,,$$
and that is a contradiction with \eq{35}. We conclude that $\{u_n\}$ converges strongly to $u$ in $E$ and Proposition 
\ref{p8} is
proved.   \endproof

\begin{lemma}\label{l3} Assume the hypotheses of Theorem \ref{t1} are fulfilled. Then
there exists $\lambda_\star>0$ such that for any $\lambda\in(0,\lambda_\star)$
there exist $\rho$, $\alpha>0$  such that $J_\lambda(u)\geq\alpha>0$
for any $u\in E$ with $\|u\|=\rho$.
\end{lemma}
\proof
By Remark 4 and conditions \eq{g} and \eq{G1} it follows that $E$ is continuously embedded in $L^{q(x)}(\Omega)$. So,
there exists a positive constant $c_1$ such that
\begin{equation}\label{32}
|u|_{q(x)}\leq c_1\cdot\|u\|,\;\;\;\forall\;u\in E\,.
\end{equation}
where by $|\cdot|_{q(x)}$ we denoted the norm on $L^{q(x)}(\Omega)$.

We fix $\rho\in(0,1)$ such that $\rho<1/c_1$. Then  relation  \eq{32} implies
$$
|u|_{q(x)}<1,\;\;\;\forall\;u\in E,\;{\rm with}\;\|u\|=\rho\,.
$$
Furthermore, relation \eq{10} applied to $\Phi(x,t)=|t|^{q(x)}$ yields
\begin{equation}\label{33}
\int_\Omega|u|^{q(x)}\;dx\leq|u|_{q(x)}^{q^-},\;\;\;\forall\;u\in E,\;{\rm with}\;\|u\|=\rho\,.
\end{equation}
Relations \eq{32} and \eq{33} imply
\begin{equation}\label{34}
\int_\Omega|u|^{q(x)}\;dx\leq c_1^{q^-}\|u\|^{q^-},\;\;\;\forall\;u\in E,\;{\rm with}\;\|u\|=\rho\,.
\end{equation}
Taking into account  relations \eq{26}, \eq{34} and \eq{G1} we deduce that
for any $u\in E$ with $\|u\|=\rho$ the following inequalities hold true
$$
J_\lambda(u)\geq\|u\|^{\phi^0}-\lambda\cdot C_2\cdot c_1^{q^-}\cdot\|u\|^{q^-}
=\rho^{q^-}(\rho^{\phi^0-q^-}-\lambda\cdot C_2\cdot c_1^{q^-})\,.
$$
By the above inequality we remark that if we define
\begin{equation}\label{34biss}
\lambda_\star=\frac{\rho^{\phi^0-q^-}}{2\cdot C_2\cdot c_1^{q^-}}\,,
\end{equation}
then for any $\lambda\in(0,\lambda_\star)$ and any $u\in E$ with $\|u\|=\rho$ there
exists $\alpha=\frac{\rho^{\phi^0}}{2}>0$ such that
$$J_\lambda(u)\geq\alpha>0\,.$$
The proof of Lemma \ref{l3} is complete.  \endproof

\begin{lemma}\label{l4} Assume the hypotheses of Theorem \ref{t1} are fulfilled. Then
there exists $\theta\in E$ such that $\theta\geq 0$,
$\theta\neq 0$ and
$J_\lambda(t\theta)<0$,
for $t>0$ small enough.
\end{lemma}
\proof
By the hypotheses of Theorem \ref{t1} we have $q^-<\phi_0$. Let $\epsilon_0>0$ be such that $q^-+\epsilon_0<\phi_0$.
On the other hand, since $q\in C(\overline\Omega)$ it follows that there exists an open set
$\Omega_0\subset\subset\Omega$ such that $|q(x)-q^-|<\epsilon_0$ for all $x\in\Omega_0$. Thus, we
conclude that $q(x)\leq q^-+\epsilon_0<\phi_0$ for all $x\in\Omega_0$.

Let $\theta\in C_0^\infty(\Omega)\subset E$ be such that ${\rm supp}(\theta)\supset\overline\Omega_0$,
$\theta(x)=1$ for all $x\in\overline\Omega_0$ and $0\leq\theta\leq 1$ in $\Omega$.

Taking into account all the above pieces of information and relations \eq{17} and \eq{G1} we have
\begin{eqnarray*}
J_\lambda(t\cdot\theta)&=&\int_{\Omega}[\Phi(x,t|\nabla\theta(x)|)+\Phi(x,t|\theta(x)|)]\;dx-\lambda
\int_\Omega G(x,t\cdot\theta(x))\;dx\\
&\leq& t^{\phi_0}\cdot\int_{\Omega}[\Phi(x,|\nabla\theta(x)|)+\Phi(x,|\theta(x)|)]\;dx-\lambda\cdot C_1\cdot
\int_\Omega t^{q(x)}|\theta|^{q(x)}\;dx\\
&\leq&t^{\phi_0}\cdot\Lambda(\theta)-\lambda\cdot C_1\cdot\int_{\Omega_0} t^{q(x)}|\theta|^{q(x)}\;dx\\
&\leq&t^{\phi_0}\cdot\Lambda(\theta)-\lambda\cdot C_1\cdot t^{q^-+\epsilon_0}\cdot\int_{\Omega_0} 
|\theta|^{q(x)}\;dx\,,
\end{eqnarray*}
for any $t\in(0,1)$, where by $|\Omega_0|$ we denoted the Lebesgue measure of $\Omega_0$.
Therefore
$$J_\lambda(t\cdot\theta)<0$$
for $t<\delta^{1/(\phi_0-q^--\epsilon_0)}$ with
$$0<\delta<\min\left\{1,\frac{\lambda\cdot C_1\cdot \int_{\Omega_0} 
|\theta|^{q(x)}\;dx}{\Lambda(\theta)}\right\}\,.$$
Finally, we point out that $\Lambda(\theta)>0$. Indeed, it is clear that
$$0<\int_{\Omega_0} |\theta|^{q(x)}\;dx\leq\int_{\Omega} |\theta|^{q(x)}\;dx\int_{\Omega} |\theta|^{q^-}\;dx
\leq c_1^{q^-}\|u\|^{q^-}\,.$$
Thus, we infer that $\|\theta\|>0$. That fact and relations \eq{25} and \eq{26} imply that $\Lambda(\theta)>0$.
The proof of Lemma \ref{l4} is complete.   \endproof
\medskip

{\sc Proof of Theorem \ref{t1}.}
Let $\lambda_\star>0$ be defined as in \eq{34biss} and $\lambda\in(0,\lambda_\star)$.
By Lemma \ref{l3} it follows that on the boundary of the  ball centered in the origin and of
radius $\rho$ in $E$, denoted by $B_\rho(0)$, we have
$$
\inf\limits_{\partial B_\rho(0)}J_\lambda>0.
$$
On the other hand, by Lemma \ref{l4}, there exists $\theta\in E$ such that
$J_\lambda(t\cdot\theta)<0$ for all $t>0$ small enough. Moreover, relations \eq{26}, \eq{34} and
\eq{G1} imply that for any $u\in B_\rho(0)$ we have
$$J_\lambda(u)\geq\|u\|^{\phi^0}-\lambda\cdot C_2\cdot c_1^{q^-}\|u\|^{q^-}\,.$$
It follows that
$$-\infty<\underline{c}:=\inf\limits_{\overline{B_\rho(0)}}J_\lambda<0\,.$$
We let now $0<\epsilon<\inf_{\partial B_\rho(0)}J_\lambda-\inf_{B_\rho(0)}J_\lambda$.
Applying Ekeland's variational principle \cite{E} to the functional
$J_\lambda:\overline{B_\rho(0)}\rightarrow\RR$,  we find
$u_\epsilon\in\overline{B_\rho(0)}$ such that
\begin{eqnarray*}
J_\lambda(u_\epsilon)&<&\inf\limits_{\overline{B_\rho(0)}}J_\lambda+\epsilon\\
J_\lambda(u_\epsilon)&<&J_\lambda(u)+\epsilon\cdot\|u-u_\epsilon\|,\;\;\;u\neq u_\epsilon.
\end{eqnarray*}
Since
$$J_\lambda(u_\epsilon)\leq\inf\limits_{\overline{B_\rho(0)}}J_\lambda+\epsilon\leq
\inf\limits_{B_\rho(0)} J_\lambda+\epsilon<\inf\limits_{\partial B_\rho(0)}J_\lambda\,,$$
we deduce that $u_\epsilon\in B_\rho(0)$. Now, we define $I_\lambda:
\overline{B_\rho(0)}\rightarrow\RR$ by $I_\lambda(u)=J_\lambda(u)+\epsilon\cdot
\|u-u_\epsilon\|$. It is clear that $u_\epsilon$ is a minimum point
of $I_\lambda$ and thus
$$\di\frac{I_\lambda(u_\epsilon+t\cdot v)-{I_\lambda}(u_\epsilon)}
{t}\geq 0$$
for small $t>0$ and any $v\in B_1(0)$. The
above relation yields
$$\di\frac{J_\lambda(u_\epsilon+t\cdot v)-J_\lambda(u_\epsilon)}{t}+
\epsilon\cdot\|v\|\geq 0.$$
Letting $t\rightarrow 0$ it follows that $\langle J_\lambda^{'}
(u_\epsilon),v\rangle+\epsilon\cdot\|v\|>0$ and we infer that
$\|J_\lambda^{'}(u_\epsilon)\|\leq\epsilon$.

We deduce that there exists a sequence
$\{w_n\}\subset B_\rho(0)$ such that
\begin{equation}\label{PS}
J_\lambda(w_n)\rightarrow{\underline c}\;\;\;{\rm and}\;\;\;
J_\lambda^{'}(w_n)\rightarrow 0.
\end{equation}
It is clear that $\{w_n\}$ is bounded in $E$. Thus, there exists
$w\in E$ such that, up  to a subsequence, $\{w_n\}$ converges
weakly to $w$ in $E$. Since, by Remark 4, $E$ is compactly
embedded in $L^{q(x)}(\Omega)$ it follows that $\{w_n\}$ converges
strongly to $w$ in $L^{q(x)}(\Omega)$. The above information
combined with relation \eq{g} and  H\"older's inequality implies
\begin{equation}\label{39}
\begin{array}{lll}
\left|\di\int_\Omega g(x,w_n)\cdot(w_n-w)\;dx\right|&\leq& C_0\cdot\di\int_\Omega |w_n|^{q(x)-1}|w_n-w|\;dx\\
&\leq& C_0\cdot|\;|w_n|^{q(x)-1}\;|_{\frac{q(x)}{q(x)-1}}\cdot|w_n-w|_{q(x)}\rightarrow 0,\;\;\;{\rm as}\;
n\rightarrow\infty\,.
\end{array}
\end{equation}
On the other hand, by \eq{PS} we have
\begin{equation}\label{40}
\lim_{n\rightarrow\infty}\langle J_\lambda^{'}(w_n),w_n-w\rangle=0\,.
\end{equation}
Relations \eq{39} and \eq{40} imply
$$\lim_{n\rightarrow\infty}\langle\Lambda^{'}(w_n),w_n-w\rangle=0\,.$$
Thus, by Proposition \ref{p8} we find that $\{w_n\}$ converges strongly to $w$ in $E$.
So, by \eq{PS},
$$
J_\lambda(w)=\underline c<0\;\;\;{\rm and}\;\;\;J_\lambda^{'}(w)=0\,.
$$
We conclude that $w$ is a nontrivial weak solution for problem \eq{1} for any
$\lambda\in(0,\lambda_\star)$. The proof of Theorem \ref{t1} is complete.  \endproof

\begin{lemma}\label{l6}
Assume the hypotheses of Theorem \ref{t2} are fulfilled. Then for any $\lambda>0$ the functional $J_\lambda$ is 
coercive.
\end{lemma}
\proof
For each $u\in E$ with $\|u\|>1$ and $\lambda>0$ relations \eq{25}, \eq{g} and Remark 4 imply
\begin{eqnarray*}
J_\lambda(u)&\geq&\|u\|^{\phi_0}-\lambda\cdot C_2\cdot\int_\Omega|u|^{q(x)}\;dx\\
&\geq&\|u\|^{\phi_0}-\lambda\cdot C_2\cdot\left[\int_\Omega|u|^{q^-}\;dx+\int_\Omega|u|^{q^+}\;dx\right]\\
&\geq&\|u\|^{\phi_0}-\lambda\cdot C_3\cdot[\|u\|^{q^-}+\|u\|^{q^+}]\,,
\end{eqnarray*}
where $C_3$ is a positive constant. Since $q^+<\phi_0$ the above inequality implies that $J_\lambda(u)\rightarrow
\infty$ as $\|u\|\rightarrow\infty$, that is, $J_\lambda$ is coercive. The proof of Lemma \ref{l6} is complete.  
\endproof
\medskip

{\sc Proof of Theorem \ref{t2}.}
Since $q^+<\phi_0$ it follows that $q^-<\phi_0$ and thus, by Theorem \ref{t1} there exists $\lambda_\star>0$ such 
that
for any $\lambda\in(0,\lambda_\star)$ problem \eq{1} has a nontrivial weak solution.

Next, by Lemma \ref{l6} and Proposition \ref{p7} we infer that $J_\lambda$ is coercive and weakly lower 
semi-continuous in $E$,
for all $\lambda>0$. Then Theorem 1.2 in \cite{S} implies that
there exists $u_\lambda\in E$ a global minimizer of $I_\lambda$ and
thus a weak solution of problem \eq{1}.

We show that $u_\lambda$ is not trivial for $\lambda$ large enough.
Indeed, letting $t_0>1$ be a fixed real and $u_0(x)=t_0$, for all $x\in\Omega$ we have $u_0\in E$ and
\begin{eqnarray*}
J_\lambda(u_0)&=&\Lambda(u_0)-\lambda\int_\Omega G(x,u_0)\;dx\\
&\leq&\int_\Omega\Phi(x,t_0)\;dx-\lambda\cdot C_1\cdot\int_\Omega|t_0|^{q(x)}\;dx\\
&\leq&L-\lambda\cdot C_1\cdot t_0^{q^+}\cdot|\Omega_1|\,,
\end{eqnarray*}
where $L$ is a positive constant.
Thus, there exists $\lambda^\star>0$ such that $J_\lambda(u_0)<0$
for any $\lambda\in[\lambda^\star,\infty)$. It follows that
$J_\lambda(u_\lambda)<0$ for any $\lambda\geq\lambda^\star$ and
thus $u_\lambda$ is a nontrivial weak solution of problem \eq{1}
for $\lambda$ large enough. The proof of Theorem \ref{t2} is
complete. \endproof

\section{Examples}
In this section we point out certain examples of functions $\phi(x,t)$ and $\Phi(x,t)$ for which the
results of this paper can be applied.

I) We can take
$$\phi(x,t)=p(x)|t|^{p(x)-2}t\;\;\;{\rm and}\;\;\;\Phi(x,t)=|t|^{p(x)}\,,$$
with $p(x)\in C(\overline\Omega)$ satisfying
$2\leq p(x)<N$, for all $x\in\overline\Omega$. It is easy to verify that $\phi$ and $\Phi$ satisfy conditions 
($\phi$),
($\Phi_1$), ($\Phi_2$), \eq{5}, \eq{7} and \eq{29} since in this case  we can take $\phi_0=p^-$ and $\phi^0=p^+$.
\smallskip

II) We can take
$$\phi(x,t)=p(x)\frac{|t|^{p(x)-2}t}{\log(1+|t|)}\;\;\;{\rm and}\;\;\; \Phi(x,t)=\frac{|t|^{p(x)}}{\log(1+|t|)}+
\int_0^{|t|}\frac{s^{p(x)}}{(1+s)(\log(1+s))^2}\;ds\,,$$
with $p(x)\in C(\overline\Omega)$ satisfying  $3\leq p(x)<N$, for all $x\in\overline\Omega$.

It is easy to see that relations ($\phi$), ($\Phi_1$) and ($\Phi_2$) are verified.

For each $x\in\overline\Omega$ fixed by Example 3 on p. 243 in \cite{Clem2} we have
$$p(x)-1\leq\frac{t\cdot\phi(x,t)}{\Phi(x,t)}\leq p(x),\;\;\;\forall\;t\geq 0\,.$$
Thus, relation \eq{5} holds true with $\phi_0=p^--1$ and $\phi^0=p^+$.

Next, $\Phi$ satisfies condition \eq{29} since
$$\Phi(x,t)\geq t^{p(x)-1},\;\;\;\forall\;x\in\overline\Omega,\;t\geq 0\,.$$

Finally, we point out that trivial computations imply that $\frac{d^2(\Phi(x,\sqrt{t}))}{dt^2}\geq 0$ for all
$x\in\overline\Omega$ and $t\geq 0$. Thus, relation \eq{7} is satisfied.
\smallskip

III) We can take
$$\phi(x,t)=p(x)\cdot\log(1+\alpha+|t|)\cdot|t|^{p(x)-1}t\,,$$
and
$$\Phi(x,t)=\log(1+\alpha+|t|)\cdot|t|^{p(x)}-\int_0^{|t|}\frac{s^{p(x)}}{1+\alpha+s}\;dx\,,$$
where $\alpha>0$ is a constant and $p(x)\in C(\overline\Omega)$ satisfying  $2\leq p(x)<N$, for all 
$x\in\overline\Omega$.

It is easy to see that relations ($\phi$), ($\Phi_1$) and ($\Phi_2$) are verified.

Next, it is easy to remark that for each $x\in\overline\Omega$ fixed we have
$$p(x)\leq\frac{t\cdot\phi(x,t)}{\Phi(x,t)},\;\;\;\forall\;t\geq 0\,.$$
The above information shows that taking $\phi_0=p^-$ we have
$$1<p^-\leq\frac{t\cdot\phi(x,t)}{\Phi(x,t)},\;\;\;\forall\;x\in\overline\Omega,\;t\geq 0\,.$$
On the other hand, some simple computations imply
$$\lim_{t\rightarrow\infty}\frac{t\cdot\phi(x,t)}{\Phi(x,t)}=p(x),\;\;\;\forall\;x\in\overline\Omega$$
and
$$\lim_{t\rightarrow 0}\frac{t\cdot\phi(x,t)}{\Phi(x,t)}=p(x),\;\;\;\forall\;x\in\overline\Omega\,.$$
Thus, defining $H(x,t)=\frac{t\cdot\phi(x,t)}{\Phi(x,t)}$ we remark that $H(x,t)$ is continuous
on $\overline\Omega\times[0,\infty)$ and $1<p^-\leq\lim_{t\rightarrow 0}H(x,t)\leq p^+<\infty$ and
$1<p^-\leq\lim_{t\rightarrow\infty}H(x,t)\leq p^+<\infty$. It follows that
$$\phi^0=\sup_{t>0,\;x\in\overline\Omega}\frac{t\cdot\phi(x,t)}{\Phi(x,t)}<\infty\,.$$
We conclude that relation \eq{5} is satisfied.

On the other hand, since
$$\phi(x,t)\geq p^-\cdot\log(1+\alpha)\cdot t^{p(x)-1},\;\;\;\forall\;x\in\overline\Omega,\;t\geq 0\,,$$
it follows that
$$\Phi(x,t)\geq\frac{p^-}{p^+}\cdot(1+\alpha)\cdot t^{p(x)},\;\;\;\forall\;x\in\overline\Omega,\;t\geq 0\,.$$
The above relation assures that relation \eq{29} is verified.

Finally, we point out that trivial computations imply that $\frac{d^2(\Phi(x,\sqrt{t}))}{dt^2}\geq 0$ for all
$x\in\overline\Omega$ and $t\geq 0$ and thus, relation \eq{7} is satisfied.

\medskip
{\bf Acknowledgements.}
The authors have been
supported by Grant CNCSIS PNII--79/2007 {\it ``Procese
Neliniare Degenerate \c si Singulare"}.

\end{document}